
\magnification=\magstep1
\documentstyle{amsppt}
\def\Vol{\operatorname{Vol}}
\def\Pos{\operatorname{Pos}}
\def\index{\operatorname{index}}
\def\coker{\operatorname{coker}}
\def\ker{\operatorname{ker}}

\def\={\overset\text{def}\to =}
\def\RP{\bR\bP}
\def\HP{\bH\bP}
\def\CP{\bC\bP}

\def\psc{positive scalar curvature}
\def\bZ{\Bbb Z}
\def\bH{\Bbb H}

\def\bP{\Bbb P}
\def\bR{\Bbb R}
\def\bC{\Bbb C}
\def\Z#1{\bZ/{#1}}

\topmatter
\leftheadtext{Jonathan Rosenberg and Stephan Stolz}
\title A ``stable" version of the Gromov-Lawson conjecture
\endtitle
\author Jonathan Rosenberg$^1$  and Stephan Stolz$^2$\endauthor
\affil $^1$Department of Mathematics,\\
University of Maryland\\
College Park, Maryland 20742, U.S.A. \\
and \\
$^2$Department of Mathematics,\\
University of Notre Dame \\
Notre Dame, Indiana 46556, U.S.A.
\endaffil
\address
Department of Mathematics, University of Maryland,
College Park, MD 20742
\endaddress
\address
Department of Mathematics, University of Notre Dame,
Notre Dame, IN 46556
\endaddress
\thanks
$^1$Partially supported by NSF Grants
DMS--90-02642 and DMS--92-25063. \endgraf
$^2$Partially supported by NSF Grant
DMS--90-02594.
\endthanks
\keywords positive scalar curvature, $KO$-theory,
bordism, Dirac operator, assembly map, real C*-algebra
\endkeywords
\subjclass Primary 53C21; Secondary 55N22,
19L41, 55N15, 58G12, 46L80, 19L64, 19B28, 57R90, 55R35, 55T15,
55U20, 19K56
\endsubjclass
\email jmr\@math.umd.edu; stolz.1\@nd.edu
\endemail
\abstract
We discuss a conjecture of Gromov and Lawson, later modified by
Rosenberg, concerning the existence of \psc\ metrics. It says that a
closed spin manifold $M$ of dimension $n\ge 5$ has a \psc\
metric if and only if the index of a suitable ``Dirac" operator
in $KO_n(C^*(\pi_1(M)))$, the real $K$@-theory of the
group $C^*$@-algebra of the fundamental group of $M$, vanishes.
It is known that the vanishing of the index is necessary for
existence of a \psc\ metric on $M$, but this is known to be
a sufficient condition only if $\pi_1(M)$ is the trivial group,
$\Z2$, an odd order cyclic group, or one of a fairly small class
of torsion-free groups.

We note that the groups $KO_n(C^*(\pi))$ are periodic in $n$ with
period $8$, whereas there is no obvious
periodicity in the original geometric problem.
This leads us to introduce a ``stable'' version of the
Gromov-Lawson conjecture, which makes the weaker statement that
the product of $M$ with enough copies of the ``Bott manifold" $B$
has a \psc\ metric if and only if the index of the Dirac operator
on $M$ vanishes. (Here $B$ is a simply connected $8$@-manifold
which represents the periodicity element in $KO_8(pt)$.)
We prove the stable Gromov-Lawson conjecture
for all spin manifolds with finite fundamental group and for many
spin manifolds with infinite fundamental group.
\endabstract
\date January, 1994\enddate
\endtopmatter
\document
\newpage
\heading Contents \endheading
\roster
\item"" Introduction
\item"1." Constructions of \psc\ metrics
\item"2." Obstructions to \psc\ metrics
\item"3." The Gromov-Lawson-Rosenberg conjecture
\item"4." The Stable Conjecture
\endroster

\head Introduction
\endhead

Let $M$ be an $n$-dimensional manifold (all manifolds considered in
this paper are
smooth, compact, and, unless otherwise specified,
they are connected and their boundary is empty).
In this note we study the following

\remark{Question} Under which (topological) conditions does $M$ admit a
\psc\ metric, i.e.\ a Riemannian
metric whose scalar curvature function is positive everywhere?
\endremark

This question has certainly a more differential geometric flavor then
one would
expect at a conference on homotopy theory. However, thanks to results
of Gromov-Lawson and (independently) Schoen-Yau, the answer to the
above question
depends only on the bordism class of $M$ in a suitable bordism group.
Via the Pontrjagin-Thom construction, this bordism group can be interpreted
as a homotopy group of a Thom spectrum, and then homotopy theoretic
techniques, notably the Adams spectral sequence, are used to get
(partial) results to the above question.

\head 1. Constructions of \psc\ metrics
\endhead

Let $g$ be a Riemannian metric on a manifold $M$ of dimension $n$.
The scalar curvature is a smooth function $s\colon M @>>> \bR$, which is
obtained from the curvature tensor by contracting twice. More geometrically,
the scalar curvature at a point
$p$ is a measure of how fast the volume of the ball of radius $r$
around $p$ is growing with $r$. More precisely, we compare
$\Vol B_r(M,p)$,
the volume of the ball of radius $r$ around $p$, with
$\Vol B_r(\bR^n,0)$, the volume of the ball of radius $r$ in
$n$-dimensional Euclidean space $\bR^n$, by expressing their quotient
as a power series in $r$. We get:
$$
\frac {\Vol B_r(M,p)}{\Vol B_r(\bR^n,0)}=1-\frac{s(p)}{6(n+2)}r^2+\ldots.
\tag 1.1$$
In particular, a Riemannian manifold has \psc\ if and only if
the volume of (small)
balls grows slower than the volume of the corresponding
Euclidean balls. Examples of manifolds with \psc\ are the
$n$-dimensional sphere $S^n$, with its usual metric,
as well as certain quotients of $S^n$, like projective
spaces (real, complex, or quaternionic),
equipped with the metric induced by the standard metric on $S^n$.
Here, we assume of course $n\ge 2$, since
the scalar curvature of a $1$@-dimensional manifold is identically zero,
as one can see immediately from equation (1.1).

One possible approach to the question of the introduction
is to study whether the answer changes when we
modify $M$. A modification of manifolds which is popular among geometric
topologists and has been very successful in the diffeomorphism
classification of manifolds is ``surgery".
Given a manifold $M$ of dimension $n$, and an embedding
$S^k\times D^{n-k}\hookrightarrow M$, we remove from $M$ the open subset
$S^k\times \dot D^{n-k}$ to get a manifold with boundary $S^k\times S^{n-k-1}$
and glue it with $D^{k+1}\times S^{n-k-1}$ along their common boundary
to get a closed manifold $\widehat M$. One says that $\widehat M$ is the
result of
a {\it $k$@-surgery} on $M$. A crucial result in the subject is the following
``Surgery Theorem" proved independently by Gromov-Lawson and Schoen-Yau.

\proclaim{1.2 Surgery Theorem
(Gromov-Lawson \cite{GL2}, Schoen-Yau \cite{SY})}

Let $M$ be an $n$-dimensional manifold (not necessarily connected)
which admits a \psc\ metric, and
suppose that $\widehat M$ is obtained from $M$ by a $k$@-surgery. If
$n-k$, the ``codimension" of the surgery, is greater than or equal to
$3$, then $\widehat M$
carries a \psc\ metric.
\endproclaim

\demo{Sketch of proof}
The key step in the proof of this result is a careful
deformation of the
original metric on $M\setminus S^k\times \dot D^{n-k}$
in a neighbourhood
of its boundary in such a way that
\roster
\item The deformed metric has still \psc\.
\item The deformed metric fits together with the ``standard metric"
on $D^{k+1}\times S^{n-k-1}$ to give a metric on $\widehat M$.
\endroster
Here the ``standard metric" is the product of the usual metrics
on $D^{k+1}$ and on $S^{n-k-1}$. We observe that
this metric has \psc\,
{\it provided\/} $n-k\ge 3$, which explains the codimension condition
in the Surgery Theorem.
\qed\enddemo

For applications of this
result, it is important to characterize those
manifolds obtainable from a given manifold $N$ by a sequence of
surgeries of codimension $\ge 3$. Morse theory implies that
a manifold $M$ can be obtained from $N$ by a sequence of surgeries
(without restrictions on the codimension)
if and only if $N$ is bordant to $M$, i.e\.
if there is a manifold $W$ whose boundary is the disjoint union of $M$ and $N$.
The crucial `if' part of this assertion can be seen as follows:
Given a
bordism $W$ between $M$ and $N$, we can find
a ``Morse function" on $W$, that is, a smooth
function $h\colon W @>>> [0,1]$ satisfying the
following conditions:
\roster
\item"(1.3a)" $h^{-1}(0)=M,\qquad h^{-1}(1)=N$.
\item"(1.3b)" The Hessian of $h$ at a critical point is non-degenerate.
\item"(1.3c)" There is
at most one critical point on
on each level set $h^{-1}(t)$, $0<t<1$,
and none for
$t=0,1$.
\endroster

Given a subinterval $[t_1,t_2]\subseteq [0,1]$
the level sets $h^{-1}(t_1)$ and
$h^{-1}(t_2)$ are diffeomorphic, provided there is no critical point in
$h^{-1}([t_1,t_2])$.
If there is exactly one critical point
in $h^{-1}([t_1,t_2])$, then
$h^{-1}(t_1)$ is obtained
from $h^{-1}(t_2)$ by a surgery whose codimension is the
index of this critical point  (cf\. \cite{Mi}). In particular, we
get the following result which we state as a lemma for future
reference.

\proclaim{1.4 Lemma} A manifold $M$ can be obtained from
another manifold $N$ by a sequence of surgeries of
codimension $\ge 3$, if we can find a bordism $W$ between them and
 a Morse function $h\colon W @>>> [0,1]$
whose critical points have index $\ge 3$.
\endproclaim

Gromov and Lawson noticed that one can find such a pair $(W,h)$,
provided $M$ and $N$ represent the same element in a suitable
bordism group. Before stating their result, we introduce some
notation. Given a topological space $X$, we denote by $\Omega_n(X)$
the bordism classes of pairs $(M^n,f)$, where $M$ is an $n$-dimensional
closed manifold, and $f\colon M @>>> X$ is a map (two such
pairs $(M_1,f_1)$, $(M_2,f_2)$
are bordant if there is a bordism $W$ between $M_1$ and $M_2$, and a map
$F\colon W @>>> X$ whose restriction to $M_i$ is $f_i$ for $i=1,2$).
Disjoint union of such pairs gives $\Omega_n(X)$ the structure of an
abelian group. We are actually interested in  variations
 of this bordism
group where all manifolds are equipped with compatible orientations resp\.
spin structures. The usual notation for these bordism groups is
$\Omega_n^{SO}(X)$ resp\. $\Omega_n^{spin}(X)$. Let $\Pos_n^{SO}(X)$
(resp\. $\Pos_n^{spin}(X)$) be the subgroup of
$\Omega_n^{SO}(X)$ (resp\. $\Omega_n^{spin}(X)$) consisting of
bordism classes represented by pairs $(M,f)$ for which $M$ admits
a \psc\ metric.

\proclaim{1.5 Bordism Theorem (Gromov-Lawson \cite{GL2})}
Let $M$ be a manifold of dimension $n\ge 5$ with fundamental group $\pi$.
Let $u\colon M @>>> B\pi$ be the classifying map of the universal
covering $\widetilde M @>>> M$. Then $M$ has a \psc\ metric if and only if
$$
[M,u]\in
\cases
\Pos_n^{spin}(B\pi)&\text{if $M$ is spin},\\
\Pos_n^{SO}(B\pi)&\text{if $M$ is oriented and $\widetilde M$ is not spin}.
\endcases
$$
\endproclaim

\remark{Remark} This result doesn't cover non-orientable manifolds or
manifolds which are non-spin, but whose universal cover is spin.
There is however a
general result \cite{RS1} covering all manifolds of dimension $n\ge 5$,
where the bordism groups $\Omega_n^{SO}(X)$
resp\. $\Omega_n^{spin}(X)$ have to be replaced
by more general bordism groups.
\endremark

\demo{Proof}
We have to show that if $(M,u)$ is bordant to $(N,f)$, and $N$ has
a \psc\ metric, then so does $M$. Combining the Surgery Theorem with
 Lemma 1.4 we see that it suffices to find a bordism $W$
between $M$ and $N$ and a Morse function $h\colon W @>>> [0,1]$
whose critical points
have index $\ge 3$.

We recall from Morse theory that $W$ is homotopy equivalent to
a space obtained by attaching cells to $M$, with each cell of dimension
$i$ corresponding to a
critical point of index $i$ of a Morse function on $W$.
In particular, if there is a
Morse function whose critical points have index $\ge 3$, then the
inclusion $M \hookrightarrow W$ is a $2$@-equivalence
(i.e\. it induces an isomorphism on the $i$@-th
homotopy group for $i<2$ and a surjection for $i=2$).
Conversely, the techniques used in the proof of the $s$@-Cobordism
Theorem show that if the dimension of $M$ is greater or equal to $5$
and the inclusion $i\colon M \hookrightarrow W$ is a $2$@-equivalence,
then one can in fact find a Morse function on $W$ whose
critical points have index $\ge 3$.

Now let's assume that $M$ is a spin manifold, and that $(W,F)$
is a bordism between $(M,u)$ and $(N,f)$.
Changing the bordism
$(W,F)$ if necessary by surgeries in the
interior of $W$, we can assume that $F\colon W @>>> B\pi$ is
a $3$@-equivalence. Since $F\circ i=u$, and $u$ is a
$2$@-equivalence, this implies that
$i$ is a $2$@-equivalence and proves the theorem in this case.

If $M$ and hence $W$ are not spin, this argument doesn't work,
since there might be a non-trivial class in $\pi_2(W)$
which can't be killed by surgery since an embedded $2$@-sphere
representing it has a non-trivial normal bundle.
In this case we replace $u$ by
$u'=u\times\nu_M\colon M @>>> B\pi \times BSO$,
where $\nu_M$ is the classifying map of the stable
normal bundle of $M$. The assumption that $\widetilde M$ is
not spin guarantees that $u'$ is a $2$@-equivalence. Moreover,
$u'=F'\circ i$, where
$F'=F\times\nu_W\colon W @>>> B\pi\times BSO$, and as before we
can make $F'$ a $3$@-equivalence by surgeries on $W$, since any
element in the kernel of $F'_*\colon \pi_2(W) @>>>
\pi_2(B\pi\times BSO)$ can be represented by an embedded $2$@-sphere
with trivial normal bundle.
\qed\enddemo

The Bordism Theorem shows that
classifying all manifolds of \psc\ of dimension $n$ with fundamental group
$\pi$, which are spin manifolds (resp\. orientable with
$\widetilde M$ non-spin) is equivalent to determining the
group $\Pos_n^{spin}(B\pi)$ (resp\. $\Pos_n^{SO}(B\pi)$).
The following two observations turn out to be very useful
for constructing elements of $\Pos_n(B\pi)$.

\proclaim{1.6 Observation}
If $M$, $N$ are manifolds, and $M$ admits a \psc\ metric, then so
does $M\times N$.
\endproclaim

\demo{Proof} Let $g$ be a \psc\ metric on $M$ and let $h$ be any
Riemannian metric on $N$. The scalar curvature of
the product metric $g\oplus h$ on $M\times N$ is not necessarily positive.
However, if we `shrink' $M$ by replacing the metric $g$ by $tg$
for some real number $0<t<1$, the scalar curvature
at a point $(p,q)$ of
$M\times N$ is given by
$$
s^{tg\oplus h}(p,q)=s^{tg}(p)+s^h(q)=\frac 1t s^g(p)+s^h(q),
$$
which is positive for $t$ sufficiently small. Since our manifolds
are compact, we can choose $t$ such that the scalar curvature of
$tg\oplus h$ is positive everywhere.
\qed\enddemo

This observation shows in particular that
$\Pos_*^{SO}(pt)\=\bigoplus_n\Pos_n^{SO}(pt)$ is an ideal in the
oriented bordism ring $\Omega_*^{SO}(pt)\=\bigoplus_n
\Omega_n^{SO}(pt)$ (multiplication induced by Cartesian product of
manifolds), and the analogous statement for $\Pos^{spin}_*(pt)$.

The ``shrinking" argument also works in the case of a ``twisted
product" \cite{St1}.

\proclaim{1.7 Proposition} Let $g$ be a \psc\ metric on a manifold $M$,
and let $p\colon E @>>> N$ be a fiber bundle with fiber $M$ whose
structure group acts on $M$ by isometries. Then $E$ admits
a \psc\ metric.
\endproclaim

Now let's consider \psc\ metrics on simply connected manifolds.
The relevant bordism group then is $\Omega_n^{spin}$ for
spin manifolds and $\Omega_n^{SO}$ for non-spin manifolds. The
oriented bordism ring $\Omega_*^{SO}$ was computed by Wall \cite{Wa},
who
also constructed manifolds which represent multiplicative
generators of this ring. The nice thing is that all these
manifolds are projective bundles of (real or complex)
vector bundles. In particular, they carry \psc\ metrics by
Proposition
1.7, and hence the Bordism Theorem implies the following result.

\proclaim{1.8 Theorem (Gromov-Lawson \cite{GL2, Cor\. C})}
Let $M$ be a simply connected manifold of dimension $n\ge 5$, which
does not admit a spin structure. Then $M$ carries a \psc\ metric.
\endproclaim

For spin manifolds the story is more complicated, since not every
spin manifold has a \psc\ metric, as we'll see in the next section.

\head 2. Obstructions to \psc\ metrics
\endhead

\proclaim{2.1 Theorem (Lichnerowicz \cite{Li}, 1963)}
Let $M$ be a spin manifold of dimension $n=4k$, which has
a \psc\ metric $g$. Then the $\hat A$-genus $\hat A(M)$ vanishes.
\endproclaim

The $\hat A$-genus of an orientable manifold
$M$ is a rational number, obtained by
evaluating a certain polynomial in the Pontrjagin classes of $M$
on the fundamental class of $M$.

\demo{Proof}
If $M$ is a spin manifold, the Atiyah-Singer Index Theorem
implies
$$
\hat A(M)=\index(D)\=\dim\,\ker\,D-\dim\,\coker\, D,
$$
where $D$ is the ``Dirac operator" on $M$.
On the other hand, it follows from the ``Weitzenb\"ock formula"
that the Dirac operator is invertible if the Riemannian metric
used in the construction of $D$ has \psc, and in particular
the index of $D$ is zero in that case.
\qed\enddemo

Later Hitchin found additional obstructions to the
existence of \psc\ metrics on spin manifolds of dimension
$n\equiv 1,2\mod{8}$ \cite{Hi}.
Hitchin uses
a generalization of the Dirac operator which is called the
``$C\ell_n$@-linear Atiyah-Singer operator"
in the book of Lawson-Michelsohn \cite{LM}. It commutes with
an action of the Clifford algebra $C\ell_n$ and has
a ``Clifford index" in $KO_n(pt)$ (cf\. \cite{LM}, Ch\. II, \S 7).
We will use the notation $\alpha(M)\in KO_n(pt)$ for the
(Clifford) index of the Atiyah-Singer operator on an $n$@-dimensional
spin manifold $M$. Making use again of the ``Weitzenb\"ock formula"
one concludes:

\proclaim{2.2 Theorem (Hitchin \cite{Hi})}
If $M$ is a spin manifold
which admits a \psc, then $\alpha(M)$ vanishes.
\endproclaim

\remark{2.3 Remark}
To compare Theorems 2.1 and 2.2, we recall that by Bott periodicity
the groups $KO_n(pt)$ are as follows.
$$
KO_n(pt)=
\cases
\bZ&n\equiv 0\mod{4}\\
\Z2&n\equiv 1,2\mod{8}\\
0&\text{otherwise}
\endcases.
$$
If $n$ is divisible by $4$, $\hat A(M)$ and $\alpha(M)$ agree
up to a factor (cf\. \cite{LM}, Ch\. II, Thm\. 7.10),
and hence Hitchin's result implies
Lichnerowicz's result. But Hitchin's result is more general, since
in dimensions $n\equiv 1,2\mod{8}$, there are spin
manifolds $M$ with $\alpha(M)\ne 0$. In fact, {\it every}
spin manifold of dimension $n\equiv 1,2\mod{8}$, $n\ge 9$,
is homeomorphic to one with non-trivial $\alpha$@-invariant
\cite{LM}, Ch\. IV, Cor\. 4.2.
This follows from the fact that the homotopy sphere $\Sigma^n$
corresponding to an element in Adams' $\mu$@-family has non-trivial
$\alpha$@-invariant, and hence for any spin manifold $M$
either $M$ or the connected sum $M\#\Sigma$ (which is homeomorphic
to $M$) has a non-trivial $\alpha$@-invariant.
This shows that the question whether a manifold admits a \psc\
metric is pretty subtle: the answer in general depends on the
differentiable structure of $M$.
\endremark

Hitchin's result was generalized by Rosenberg who constructs a
``Dirac operator" whose index is an element of the $K$@-theory
of a $C^*$@-algebra. More precisely, let $M$ be a spin manifold
of dimension $n$, and let $f\colon M @>>> B\pi$ be a map to
the classifying space of a discrete group $\pi$
(not necessarily the
fundamental group of $M$, but that is the main case of interest).
Then Rosenberg constructs a ``Dirac operator"
$D$ with $\index(D)$ in $KO_n(C^*(\pi))\= \pi_n(BGL(C^*(\pi)))$.
Here $C^*(\pi)$ is the (reduced) group $C^*$@-algebra of $\pi$,
which is a suitable completion of the real group ring $\bR \pi$,
and $BGL(C^*(\pi))$ is the classifying space of the general linear
group of $C^*(\pi)$ (regarded as a topological group).

It turns out that $\index(D)\in KO_n(C^*(\pi))$
is independent of the metric used in the construction of $D$.
We will use the notation $\alpha(M,f)$ for $\index (D)$, which
agrees with $\alpha(M)$ if $\pi$ is the trivial group.
The Weitzenb\"ock formula shows again
 that the index vanishes if the metric has \psc, which
proves:

\proclaim{2.4 Theorem (Rosenberg \cite{Ro3})}
If $M$ is a spin manifold
of \psc, and $f\colon
M @>>> B\pi$ is a map to the classifying space of a discrete
group $\pi$, then $\alpha(M,f)$ vanishes.
\endproclaim

\head 3. The Gromov-Lawson-Rosenberg conjecture
\endhead

\proclaim{3.1 Conjecture (Gromov-Lawson \cite{GL3},
Rosenberg \cite{Ro2})}
Let $M$ be a spin manifold of dimension $n\ge 5$ with fundamental
group $\pi$, and let
$u\colon M @>>> B\pi$ be the classifying map of the universal covering
$\widetilde M @>>> M$. Then $M$ has a \psc\ metric if and only if
the element $\alpha(M,u)$ in $KO_n(C^*(\pi))$ vanishes.
\endproclaim

Using the Bordism Theorem 1.5 and the fact that $\alpha(M,f)$ depends
only on the bordism class $[M,f]$ we can formulate the conjecture
equivalently as follows:

\proclaim{3.2 Conjecture} $\Pos_n^{spin}(B\pi)$ is the
kernel of the homomorphism
$$
\alpha\colon \Omega_n^{spin}(B\pi)@>>> KO_n(C^*(\pi)).
$$
\endproclaim

We note that $\Pos_n^{spin}(B\pi)$ is contained in $\ker\alpha$ by
Theorem 2.4. To prove the converse, one has to show that every
class in $\ker\alpha$ can be represented by a manifold with
\psc.
The problem is that even for $\pi=\{1\}$ we don't have explicit
manifolds which generate $\ker\alpha$ (or $\Omega^{spin}_*$, for that
matter, unlike the case of the oriented bordism ring).
There are however partial results:
the kernel of $\alpha$ is trivial in
dimensions $n<8$, and is the infinite cyclic group generated by
the bordism class of
the quaternionic projective space $\HP^2$ for $n=8$.
It follows that the conjecture is true for $n\le 8$.
By finding explicit \psc\ manifolds generating $\ker\alpha$ for $n<23$,
  Rosenberg \cite{Ro3} proved the conjecture in that range.
Similarly,
Miyazaki \cite{Miy} found \psc\ manifolds
generating $\ker\alpha\otimes\bZ [\frac 12]$, thus proving the
conjecture after ``inverting $2$".

Later Stolz gave the following
characterization of the kernel of $\alpha$.

\proclaim{3.3 Theorem (Stolz \cite{St1})}
The kernel of $\alpha\colon \Omega_n^{spin} @>>> KO_n(pt)$ is equal
to the subgroup $T_n$ consisting of those bordism classes
represented by total spaces of $\HP^2$@-bundles, i.e\. fiber bundles
with fiber $\HP^2$ and structure group the projective
symplectic group $PSp(3)$.
\endproclaim

We note that the group $G=PSp(3)$ acts by isometries on $\HP^2$
equipped with its standard metric. Hence
these total spaces have \psc\ metrics by Proposition 1.7, and
thus the above result
 implies the conjecture in the simply connected case.
In \cite{St1}, the equality $[\ker\alpha]_n=T_n$ is actually only
proved localized at $2$. This is enough for the proof
of the conjecture by Miyazaki's result. Localized at odd primes,
the above theorem was proved in \cite{KrSt}, using explicit
manifold constructions ($\HP^2$@-bundles over products of
two quaternionic projective spaces).

\demo{Sketch of proof}
The proof at the prime $2$ is much more indirect and proceeds by
first translating the problem into stable homotopy theory.

The group $T_n$ can be described as the image of the homomorphism
$$
\Psi\colon \Omega_{n-8}^{spin}(BG) @>>> \Omega_n^{spin},
$$
which maps a bordism class $[N^{n-8},f]$ to $[\widehat N]$,
where $\widehat N @>>> N$ is the pull back of the `universal'
$\HP^2$@-bundle $EG\times_G\HP^2 @>>> BG$ via $f$.

Now consider the following commutative diagram.
$$
\CD
\Omega_{n-8}^{spin}(BG) @>\Psi>> \Omega_n^{spin} @>\alpha>>KO_n(pt)\\
@V\cong VV @V\cong VV @V\cong VV\\
\pi_n(MSpin\wedge\Sigma^8BG_+) @>T_*>>\pi_n(MSpin) @>D_*>>\pi_n(ko)\\
@V\widehat T_* VV @| @|\\
\pi_n(\widehat{MSpin}) @>>>\pi_n(MSpin) @>D_*>>\pi_n(ko).
\endCD
\tag 3.4$$
Here the second row is the homotopy theoretic interpretation
of the first row; the left and middle vertical
isomorphism is given by the Pontrjagin-Thom construction.
We recall that the Pontrjagin-Thom construction
gives an isomorphism between the spin bordism
group $\Omega_n^{spin}(X)$ of a space $X$ and the $n$-th
homotopy group of the spectrum
$MSpin\wedge X_+$, where $MSpin$ is the Thom spectrum
over $BSpin$, and $X_+$ is the union of $X$ and a disjoint base
point. For $n\ge 0$,
$KO_n(pt)$ is isomorphic to the $n$-th homotopy group
of a connective spectrum $ko$,
called the {\it connective real $K$@-theory
spectrum}. Moreover, there are spectrum maps $T$ (a kind of
transfer map) and $D$ (the $ko$@-orientation of $MSpin$),
whose induced maps in homotopy make the upper two squares
commutative.

 Using the families index theorem,
it is proved that the composition $D\circ T$ is homotopic to the
constant map, and hence $T$ factors through a map
$\widehat T$ from $MSpin\wedge\Sigma^8BG_+$ to
$\widehat {MSpin}$, the homotopy fiber of $D$. The
bottom row is part of the long exact homotopy sequence of this
 fibration.

The statement of the theorem is equivalent to saying
that the top row of the diagram is exact. By
exactness of the bottom row it
suffices to show that
$\widehat T_*$ is surjective at the prime $2$. The proof of this
uses the following facts.
\roster
\item The homomorphism induced by $\widehat T$ on
$\Z2$@-cohomology is a split injection of modules over the
Steenrod algebra.
\item The Adams spectral sequence converging to the
$2$@-local homotopy groups of $MSpin\wedge BG_+$ collapses.
\endroster
By (1) the map of Adams spectral sequences
induced by $\widehat T$ is a surjection on the $E_2$@-level, by (2)
on the $E_\infty$@-level, and hence $\widehat T$
is surjective on  the $2$@-local homotopy
groups.
\qed\enddemo

Now we turn to discuss Conjecture 3.1 in the case of a
non-trivial fundamental group. It might be tempting to think that
a manifold $M$ has a \psc\ metric if and only if its universal covering
does, at least if $\pi_1(M)$ is finite, by arguing one should be
able to ``average" a \psc\ metric on the universal covering of $M$ to
get a $\pi_1(M)$@-equivariant \psc\ metric which then descends
to a \psc\ metric on $M$. However, averaging a \psc\ metric
might not give a \psc\ metric as the
example below shows. On the other
hand, if $\pi_1(M)$ is finite of {\it odd order}, then it can
be shown that the vanishing of  $\alpha(\widetilde M)$ implies the
vanishing of $\alpha(M,u)$, and hence Conjecture 3.1 claims
in this case that $M$ has a \psc\ metric if and only if
$\widetilde M$ does!

\remark{3.5 Example (Bergery-Berard \cite{BB}, Example 9.1)}
Let $\Sigma$ be a $9$@-di\-men\-sional
homotopy sphere with $\alpha(\Sigma)\ne 0$ (cf\. Remark 2.3), and
let $M$ be the connected sum $(\RP^7\times S^2)\#\Sigma$.
We note that the real projective space $\RP^7$ and hence
$\RP^7\times S^2$ are spin manifolds. It follows that the
$\alpha$@-invariant of $\RP^7\times S^2$ vanishes, since it
 has a
\psc\ metric. Noting
that $\alpha$ is additive with respect to connected sum ($\alpha$
is a bordism invariant, and the connected sum is bordant to the
disjoint union), we get:
$$
\alpha(M)
=\alpha(\RP^7\times S^2)+\alpha(\Sigma)\ne 0.
$$
 Hence by Theorem 2.2 $M$ does not admit a \psc\
metric. On the other hand, we have
$$
\widetilde M\cong S^7\times S^2\#\Sigma\#\Sigma\cong S^7\times S^2,
$$
since the group of $9$@-dimensional homotopy spheres is isomorphic to
$(\Z2)^3$, and hence
 $\Sigma\#\Sigma$ is diffeomorphic to $S^7$. So we see that
$\widetilde M$ has a \psc\ metric.
\endremark

Attempting to prove Conjecture 3.1/3.2, an understanding of
$\ker \alpha$ is crucial. Hence the following factorization of
$\alpha$ is useful.
$$
\Omega_n^{spin}(B\pi) @>D_*>> ko_n(B\pi) @>p>> KO_n(B\pi)
@>A>> KO_n(C^*(\pi))
\tag 3.6$$
Here $p$ is the canonical map from connective to periodic
$KO$@-homology, and $A$ is the ``assembly map".

\proclaim{3.7 Theorem (Jung \cite{Ju}, Stolz \cite{St2})}
Let $M$ be a spin manifold of dimension $n\ge 5$ with fundamental group $\pi$.
Let $u\colon M @>>> B\pi$ be the classifying map of the universal
covering $\widetilde M @>>> M$. Then $M$ has a \psc\ metric if and only if
$D_*([M,u])\in \Pos_n^{ko}(B\pi)$, where $\Pos_n^{ko}(B\pi)$ is the image
of $D_*$ restricted to $\Pos_n^{spin}(B\pi)$.
\endproclaim

\demo{Sketch of proof}
It suffices to show $\ker D_*\subseteq \Pos_n^{spin}(B\pi)$. Away from the
prime $2$ this is proved by Jung, who gives a Baas-Sullivan description of
$ko_*(X)\otimes \bZ[\frac 12]$. In particular, $[M,u]\in\ker D_*$ implies
that the connected sum of $2^r$ copies of $(M,u)$ for some $r$
bounds a manifold with
singularities, and Jung uses this to construct a \psc\ metric on the
connected sum.

The result at the prime $2$ is due to Stolz, who proves that an odd
multiple of a bordism class  in the kernel of $D_*$ can be represented
by the total space of an $\HP^2$@-bundle. This boils down to the
homotopy theoretic statement that the middle row in diagram 3.4 is still
exact after smashing with $B\pi_+$. This follows from the fact that the
map $\widehat T$ is a split surjection of spectra, which in turn is
proved by splitting the spectra $MSpin\wedge\Sigma^8BG_+$ and
$\widehat {MSpin}$ using Adams spectral sequence techniques.
\qed\enddemo

The above result is a significant improvement compared
to the Bordism Theorem 1.5, since the connective $KO$@-theory groups
$ko_*(B\pi)$ are {\it much} smaller than $\Omega^{spin}_*(B\pi)$.
Unfortunately, $ko_*(B\pi)$ has been computed for only a handful
of finite groups, notably cyclic groups \cite{Ha}, elementary
abelian $2$-groups \cite{Yu},
and the quaternion and dihedral group of order $8$ \cite{Ba}.
Still, we immediately obtain the following.

\proclaim{3.8 Corollary}
The conjecture 3.1 is true if $p$ and $A$ are injective.
\endproclaim

Note that many torsion-free groups $\pi$ for which $A$ is injective
were listed in \cite{Ro3}. (This is related to the Novikov
conjecture, as we will note below.) If in addition $B\pi$ is stably a
wedge of spheres, then $p$ is clearly a split injection. Thus
Corollary 3.8 applies to free groups, free abelian groups,
fundamental groups of orientable surfaces, and many similar
examples.

\proclaim{3.9 Theorem (Rosenberg-Stolz \cite{RS1}, Thm\. 5.3(4))}
The conjecture 3.1 is true for $\pi\cong \Z2$.
\endproclaim

\demo{Sketch of proof}
An Adams spectral sequence calculation shows that the kernel of
$A\circ p\colon ko_n(B\pi) @>>> KO_n(C^*(\pi))$ is trivial for
$n\not\equiv 3\mod{4}$, and a finite cyclic group generated by
$D_*[\RP^n,u]$ for $n\equiv 3\mod{4}$. This implies the conjecture
by Theorem 3.7.
\qed\enddemo

\head 4. The stable conjecture
\endhead

The real $K$@-theory groups of a real $C^*$@-algebra $A$ are $8$@-periodic.
Moreover, the isomorphism $KO_n(A) @>\cong>> KO_{n+8}(A)$
is given by multiplication by the generator $b$ of
$KO_8(pt)=KO_8(\bR)\cong \bZ$.
We can find a simply connected spin manifold $B$ of dimension $8$
with $\alpha(B)=b$. There are many possible choices for $B$, but
we just pick one, and call it the ``Bott manifold".
The Bott periodicity for $KO_*(C^*(\pi))$ shows that given a manifold $M$
with fundamental group $\pi$, $\alpha(M,u)$ vanishes if and only if
$\alpha(M\times B,u)$ vanishes (we use the letter $u$ for the
classifying map of the universal covering of whatever manifold
we are talking about). This shows that Conjecture 3.1 is
equivalent to the following two conjectures:

\proclaim{4.1 Cancellation Conjecture}
Let $M$ be a spin manifold of dimension $n\ge 5$. Then $M$ has a \psc\ metric
if and only if $M\times B$ does.
\endproclaim

\proclaim{4.2 Stable Conjecture}
Let $M$ be a spin manifold.
Then $M$ has stably a \psc\ metric (i.e\. the product of $M$ with
sufficiently many copies of $B$ has a
\psc\ metric) if and only if $\alpha(M,u)=0$.
\endproclaim

An important tool for proving the Stable Conjecture is the
following geometric description of periodic $KO$@-homology.

\proclaim{4.3 Theorem (Kreck-Stolz \cite{KrSt}, Thm\. C)}
Given a space $X$, let $T_*(X)$ be the subgroup of $\Omega_*^{spin}(X)$
represented by pairs $(\widehat N,f\circ p)$, where $\widehat N @>p>> N$
is an $\HP^2$@-bundle, and $f\colon N @>>> X$ is a map. Let
$b\in \Omega_8^{spin}(pt)/T_8(pt)\cong ko_8(pt)\cong \bZ$ be the
generator, i.e., the class represented by
the Bott manifold. Then the homomorphism
$p\circ D_*\colon \Omega_*^{spin}(X) @>>> KO_*(X)$ induces an isomorphism
between
the groups $\Omega_*^{spin}(X)/T_*(X)[b^{-1}]$
and $KO_*(X)$.
\endproclaim

This shows in particular that if $[M,u]\in \Omega_*^{spin}(X)$ is in the
kernel of $p\circ D_*$, then the product $M\times B\times\ldots\times B$
with sufficiently many copies of $B$ represents an element in $T_n(X)$,
and hence carries a \psc\ metric. This implies the following result.

\proclaim{4.4 Corollary} If the assembly map
$A\colon KO_n(B\pi) @>>> KO_n(C^*(\pi))$ is injective, then the
Stable Conjecture 4.2 is true.
\endproclaim

The Novikov conjecture (or rather, a form of it) claims that $A$ is injective
for torsion free groups. It has been proved for many groups, notably
for torsion free, discrete subgroups of Lie groups \cite{Ka}.
The assembly map is definitely not injective for some groups,
e.g., finite groups. In general, we obtain the following consequence of
Theorem 4.3.

\proclaim{4.5 Corollary}
Let $M$ be a manifold of dimension $n\ge 5$ with fundamental group $\pi$.
Let $u\colon M @>>> B\pi$ be the classifying map of the universal
covering $\widetilde M @>>> M$. Then $M$ has stably
a \psc\ metric if and only if
$p\circ D_*([M,u])$ is in $\Pos_n^{KO}(B\pi)$,
where $\Pos_n^{KO}(B\pi)\subset KO_n(B\pi)$ is the image
of $\bigoplus_{k\ge 0}\Pos_{n+8k}^{spin}(B\pi)$ under  $p\circ D_*$
(here we identify $KO_{n+8k}(X)$ with $KO_n(X)$ using periodicity).
\endproclaim

We would like to point out the formal similarities between this result,
the Bordism Theorem 1.5, and Theorem 3.7. However,
$KO_*(B\pi)$ is much easier to compute than $ko_*(B\pi)$ or
$\Omega_*^{spin}(B\pi)$. For example, $KO_*(B\pi)$ for a finite group
$\pi$ can be expressed in terms of the representation ring of
$\pi$, in a fashion
similar to Atiyah-Segal's calculation of $KO^*(B\pi)$.
In fact, we obtain the description of $KO_*(B\pi)$ by
dualizing their result.

\proclaim{4.6 Theorem (Rosenberg-Stolz \cite{RS2})}
The Stable Conjecture 4.2 is true for finite groups $\pi$.
\endproclaim

\demo{Sketch of proof}
By Corollary 4.5 it suffices to show
$\ker A\subseteq \Pos_*^{KO}(B\pi)$.
Using a result of Kwasik-Schultz \cite{KwSc}, (cf\. \cite{RS1}, Prop. 5.2)
 we can assume that $\pi$ is a $p$@-group.
Also it is enough to work with $\widetilde{KO}_*(B\pi)$, since
$\ker A=\Pos_*^{KO}(B\pi)$ is true if $\pi$ is the trivial group.
 As mentioned before, the groups $\widetilde{KO}_*(B\pi)$ for
a $p$@-group $\pi$ can be calculated. The result says in particular
that
$\widetilde{KO}_*(B\pi)$ is a direct sum of finitely many copies of
$\Z{p^\infty}$ plus, for $p=2$, finitely many copies of $\Z2$.
Moreover, the kernel of the assembly map restricted to
$\widetilde{KO}_*(B\pi)$ consists precisely of the
$\Z{p^\infty}$'s, and we have to show that all those elements
can be represented by \psc\ manifolds.

This is proved first for cyclic groups $\pi$ of order $k=p^r$.
We note that the classifying space $B\Z k$ can be identified with
the sphere bundle $S(H^{\otimes k})$ of the $k$@-th tensor power
of the Hopf line bundle over complex projective space $\CP^\infty$.
Then the homotopy cofibration
$$
B\Z k=S(H^{\otimes k}) @>p>> \CP^\infty @>>> T(H^{\otimes k}),
$$
where $T(H^{\otimes k})$ denotes the Thom space of $H^{\otimes k}$,
induces a long exact sequence of $KO$@-homology groups
$$
 @>>> KO_{n+1}(T(H^{\otimes k})) @>\partial>>
\widetilde{KO}_n(B\Z k) @>p_*>>
\widetilde{KO}_n(\CP^\infty) @>>>.
$$
The group $\widetilde{KO}_n(B\Z k)$ is a torsion group, while
$\widetilde{KO}_n(\CP^\infty)$ is torsion free. Hence $p_*$
is trivial, and $\partial$ is surjective.

Using again the geometric description of $KO_*$@-homology in
Theorem 4.3, one can show that the image of $\partial$ can be
represented by total spaces of $S^1$@-bundles over simply
connected manifolds. Moreover, these manifolds are non-spin
if $p$ is odd. Hence in this case these manifolds and consequently
 the total spaces of the $S^1$@-bundles over them
admit \psc\ metrics,
which shows $\widetilde{KO}_n(B\Z k)=\widetilde{\Pos}_n^{KO}(B\Z k)$.

For $p=2$, $H^{\otimes k}$ is a spin bundle and hence using the Thom
isomorphism we get
$$
KO_{n+1}(T(H^{\otimes k}))\cong KO_{n-1}(\CP^\infty)
\cong \widetilde{KO}_n(\CP^\infty)\oplus KO_n(pt).
$$
The image of $\partial$ restricted to $KO_n(pt)$ is a finite
subgroup, and hence $\partial$ restricted to
$\widetilde{KO}_n(\CP^\infty)$ is still surjective on the
$\Z{p^\infty}$ summands. But the image of $\partial$ restricted to
$\widetilde{KO}_n(\CP^\infty)$ is again represented by total
spaces of $S^1$@-bundles over manifolds with \psc.
This proves the theorem in the case of cyclic groups.

For the general case we use the fact that for a finite group
$\pi$ the representations induced up from cyclic subgroups
span a subgroup of finite index of the representation ring of
$\pi$. It follows that the image of $\bigoplus_H KO_*(BH)$ in
$KO_*(B\pi)$, where $H$ runs through all cyclic subgroups of $\pi$,
has finite index. In particular the kernel of $\alpha$ is in the
image, which reduces the general case to the cyclic case.
\qed\enddemo

So far, we only discussed what is known about {\it spin}
manifolds with \psc\ and the group $\Pos_n^{spin}(B\pi)$, but not
the analogous group $\Pos_n^{SO}(B\pi)$, which according to
the Bordism Theorem 1.5 determines whether an
orientable manifold, whose universal cover is non-spin, has a \psc\
metric. We observe that for such a manifold $M$ the product
$M\times B$ always has a \psc\ metric, since $B$ represents the
same bordism class in $\Omega_8^{SO}(pt)$ as the disjoint union
of $64$ copies of $\CP^2\times\CP^2$. This implies
$$
[M\times B]=[M\times\coprod_{64}\CP^2\times\CP^2]\in \Pos^{SO}_{n+8}(B\pi),
$$
since $\coprod_{64}\CP^2\times\CP^2$ and hence
$M\times\coprod_{64}\CP^2\times\CP^2$ has a \psc\ metric. Then the
Bordism Theorem shows that $M\times B$ has a \psc\ metric.

This shows in particular, that methods using suitable Dirac operators and
their indices in some $KO$@-groups don't work here to give obstructions
for \psc\ metrics on $M$. One might object here that $M$ is
not spin, and hence there is no Dirac operator on $M$ anyway. However,
one can define a Dirac operator \cite{St3}, which
generalizes Rosenberg's Dirac operator and
whose index lives in the $KO$@-theory of a ``twisted" version
of $C^*(\pi_1(M))$,
were the ``twist" is determined by the first two Stiefel-Whitney classes
of $M$. The manifold can be non-spin, or even non-orientable;
the only condition needed for the construction of this operator is that
the universal covering $\widetilde M$ is spin.

At this point, one might believe that every orientable manifold $M$
of dimension $n\ge 5$, whose universal cover is non-spin, has
in fact a \psc\ metric, since this is true in the simply connected
case by Theorem 1.8, and there are no obstructions coming from the indices
of Dirac operators. However, there is another technique for producing
obstructions to \psc\ metrics, namely the minimal hypersurface
method of Schoen-Yau \cite{SY}. Using it one can show that the
connected sum $M=T^6\#(\CP^2\times S^2)$ of the $6$@-dimensional
torus with $\CP^2\times S^2$ does not admit a \psc\ metric
\cite{GL3, p.\ 186},
despite that fact that the universal covering $\widetilde M$ is non-spin.
At this point, the characterization of the subgroup
$\Pos_n^{SO}(B\pi)\subseteq \Omega_n^{SO}(B\pi)$ is wide open ---
to the authors' best knowledge there is not even a conjecture about it.

\enddocument